\documentclass[sigconf]{acmart}

\usepackage[utf8]{inputenc}

\usepackage{amsmath,amssymb}
\usepackage{amsthm}
\usepackage{mathrsfs}
\usepackage{stmaryrd}
\usepackage{enumerate}
\usepackage{comment}
\usepackage{multirow}
\usepackage{xspace}
\usepackage{tabularx,multicol}
\usepackage{tikz}
\usetikzlibrary{patterns}

\definecolor{input}{HTML}{303060}
\definecolor{output}{HTML}{804000}
\definecolor{string}{HTML}{A02020}
\definecolor{ring}{HTML}{A020A0}
\definecolor{function}{HTML}{205080} 
\definecolor{constructor}{HTML}{205080}
\definecolor{method}{HTML}{205080}
\definecolor{keyword}{HTML}{008000}
\definecolor{error}{HTML}{B01010}
\definecolor{comment}{HTML}{606060}

\newcommand{\noopsort}[1]{}

\newcommand{\ph}{\vphantom{$A^A_A$}}

\newcommand{\Z}{\mathbb Z}
\newcommand{\Zp}{\Z_p}
\newcommand{\Q}{\mathbb Q}
\newcommand{\Qp}{\Q_p}

\newcommand{\calV}{\mathcal{V}}
\newcommand{\ttv}{\texttt{v}\xspace}
\newcommand{\ttw}{\texttt{w}\xspace}

\newcommand{\calU}{\mathcal{U}}
\newcommand{\calK}{\mathcal{K}}

\newcommand{\sage}{\textsc{SageMath}\xspace}
\newcommand{\ZpCR}{\text{\color{output} \rm \tt ZpCR}\xspace}

\newcommand{\ZpL}{\text{\color{output} \rm \tt ZpL}\xspace}
\newcommand{\ZpLC}{\text{\color{output} \rm \tt ZpLC}\xspace}
\newcommand{\ZpLF}{\text{\color{output} \rm \tt ZpLF}\xspace}
\newcommand{\ZpXX}{\text{\color{output} \rm \tt ZpXX}\xspace}

\newcommand{\ind}{\text{ind}}
\newcommand{\coind}{\text{coind}}
\newcommand{\hgt}{\text{hgt}}
\newcommand{\cohgt}{\text{cohgt}}
\newcommand{\inp}{\text{in}}
\newcommand{\out}{\text{out}}

\newcommand{\cIn}{{\color{input} \tt \phantom{Zp}In}:}

\newcommand{\cZpCR}{{\color{output} \tt ZpCR}:}

\newcommand{\cZpLC}{{\color{output} \tt ZpLC}:}

\definecolor{purple}{rgb}{0.6,0,0.6}

\definecolor{answer}{rgb}{0,0.5,0.2}

\newtheorem{theo}{Theorem}[section]

\newtheorem{prop}[theo]{Proposition}

\theoremstyle{definition}
\newtheorem{rmk}[theo]{Remark}

\newtheorem{deftn}[theo]{Definition}

\begin{document}

\title{\texorpdfstring{\ZpL}{ZpL}: a p-adic precision package}

\author{Xavier Caruso}
  \affiliation{Universit\'e Rennes 1; \\
}
\email{xavier.caruso@normalesup.org}
\author{David Roe}
  \affiliation{MIT; \\
}
\email{roed@mit.edu}
\author{Tristan Vaccon}
  \affiliation{Universit\'e de Limoges; \\
  }
\email{tristan.vaccon@unilim.fr}



\begin{abstract}
We present a new package \ZpL for the mathematical software system \sage. It 
implements a sharp tracking of precision on $p$-adic numbers, following 
the theory of ultrametric precision introduced in \cite{caruso-roe-vaccon:14a}. The 
underlying algorithms are mostly based on automatic differentiation 
techniques. We introduce them, study their complexity and discuss our 
design choices.
We illustrate the benefits of our package (in comparison with previous 
implementations) with a large sample of examples coming from linear 
algebra, commutative algebra and differential equations.
\end{abstract}

\maketitle

\section{Introduction}

When computing with real and $p$-adic fields, exact results are usually impossible,
since most elements have infinite decimal or $p$-adic expansions.
Working with these fields thus requires an analysis of how precision evolves
through the sequence of steps involved in carrying out a computation.
In this paper, we describe a package for computing with $p$-adic rings and fields,
based on a series of papers by the same authors
\cite{caruso-roe-vaccon:14a,caruso-roe-vaccon:15,caruso-roe-vaccon:16,caruso-roe-vaccon:17}.
The core of the package is a method for tracking precision using $p$-adic lattices which
can yield dramatically more precise results, at the cost of increased runtime and memory usage.

The standard method for handling precision when computing with real numbers
is floating point arithmetic, which may also be used in $p$-adic computation.
At a given precision level, a finite set of representable numbers are chosen,
and arithmetic operations are defined to give a representable number
that is close to the true result \cite{ieee754:2008}.  Floating point arithmetic
has the benefit of efficient arithmetic operations, but users are responsible
for tracking the precision of the results.  Numerically unstable algorithms can
lead to very inaccurate answers \cite{higham:2002}.

If provably correct results are desired, interval arithmetic provides
an alternative to floating point.  Instead of just tracking an approximation to the answer,
the package also tracks a radius within which the true result lies.
This method is commonly used for $p$-adic computations since
the ultrametric property of $p$-adic fields frequently keeps the radius small.
Computations remain fairly efficient with this approach, but numerical instability
can still lead to dramatic losses in precision (see \S \ref{sec:demo} for many examples).
Tracking the precision of multiple variables concurrently, the set of possible
true values associated to an inexact value takes the form of an ellipsoid
with axes parallel to the coordinate axes.

For better control of precision, we may allow arbitrary axes.  This 
change would have little utility for real numbers, since such ellipsoids 
are not preserved by most functions. For $p$-adic fields, in contrast, 
differentiable maps with surjective differential will send sufficiently 
small ellipsoids to other ellipsoids.  From an algebraic perspective, 
these ellipsoids are just cosets of a lattice $H$ inside a $p$-adic 
vector space, and the main result of \cite{caruso-roe-vaccon:14a} (see 
also Proposition \ref{prop:precision} below) describes how the image of 
such a coset under a map $f$ is given exactly by applying the 
differential of $f$ to $H$.

In this paper, we describe an implementation of this idea in \sage \cite{sage}.
Rather than attaching a precision to each element, we store the precision of
many elements together by tracking a precision module for the whole collection
of variables.  As variables are created and destroyed, we update a matrix whose
rows represent the vectors in the module.  Information about the precision
of elements is extracted from the matrix as necessary.

The article is structured as follows. In \S \ref{sec:demo} we provide a demonstration
of the package, showing how it can provide more precise answers than the
traditional methods for tracking $p$-adic precision.  In particular, \S \ref{ssec:elem}
describes elementary arithmetic and the SOMOS-4 sequence, \S \ref{ssec:linalg}
gives examples from linear algebra, \S \ref{ssec:comalg} examples using polynomials,
and \S \ref{ssec:diffeq} examples of differential equations.

In \S \ref{sec:implementation} we give more details on the implementation.
\S \ref{ssec:preclemma} contains a brief overview on the theory of
$p$-adic precision of \cite{caruso-roe-vaccon:14a}. 
In the next two subsections, we explain in more details how \ZpLC and 
\ZpLF work. \S \ref{ssec:trackprec} is devoted to the implementation of 
automatic differentiation leading to the actual computation of the
module that models the precision.
In \S \ref{ssec:viewprec}, we explain how precision on any individual number
can be recovered and discuss the validity
of our results.
The complexity overhead induced by our package is analyzed in
\S \ref{ssec:complexity}.

Finally, \S \ref{sec:conclusion} contains a discussion of how we see this
package fitting into the existing $p$-adic implementations.  While these
methods do introduce overhead, they are well suited to exploring
precision behavior when designing algorithms, and can provide
hints as to when further precision analysis would be useful.

\section{Short demonstration}
\label{sec:demo}

The first step is to define the parents: the rings of $p$-adic
numbers we will work with.

\smallskip

{\noindent \small
\begin{tabular}{rl}
\cIn
 & {\color{ring}\verb?Z2?}\verb? = ?{\color{output}\verb?ZpXX?}\verb?(2, print_mode=?{\color{string}\verb?'digits'?}\verb?)? \\
 & {\color{ring}\verb?Q2?}\verb? = ?{\color{output}\verb?QpXX?}\verb?(2, print_mode=?{\color{string}\verb?'digits'?}\verb?)? \\
\end{tabular}}

\smallskip

\noindent
\ZpXX is a generic notation for \ZpCR, \ZpLC and \ZpLF.
The first, \ZpCR, is the usual constructor for $p$-adic parents in 
\sage. It tracks precision using interval arithmetic. On the contrary
\ZpLC and \ZpLF are provided by our package. In the sequel, we will
compare the outputs provided by each parent.
Results for \ZpLF are only displayed when they differ from \ZpLC.

\subsection{Elementary arithmetic}
\label{ssec:elem}

We begin our tour of the features of the \ZpL package with some basic 
arithmetic computations.
We first pick some random element $x$.
The function {\color{function}\verb?random_element?} is designed so
that it guarantees that the picked random element is the same for
each constructor \ZpCR, \ZpLC and \ZpLF.

\smallskip

{\noindent \small
\begin{tabular}{rl}
\cIn
 & \verb?x = ?{\color{function}\verb?random_element?}\verb?(?{\color{ring}\verb?Z3?}\verb?, prec=5); x? \\
\cZpCR
 & \verb?...11111? \\
\cZpLC
 & \verb?...11111? \\
\end{tabular}}

\smallskip

\noindent
Multiplication by $p$ (here $3$) is a shift on the digits and thus
leads to a gain of one digit in absolute precision.
In the example below, we observe that when this multiplication is
split into several steps, \ZpCR does not see the gain of precision
while \ZpL does.

\smallskip

{\noindent \small
\begin{tabular}{rl@{\hspace{1.5cm}}rl}
\cIn
 & \verb?3*x? &
\cIn
 & \verb?x + x + x? \\
\cZpCR
 & \verb?...111110? &
\cZpCR
 & \verb? ...11110? \\
\cZpLC
 & \verb?...111110? &
\cZpLC
 & \verb?...111110? \\
\end{tabular}}

\smallskip

%
%
%

\noindent
The same phenomenon occurs for multiplication.

\smallskip

{\noindent \small
\begin{tabular}{rl@{\hspace{1.5cm}}rl}
\cIn
 & \verb?x^3? &
\cIn
 & \verb?x * x * x? \\
\cZpCR
 & \verb?...010101? &
\cZpCR
 & \verb? ...10101? \\
\cZpLC
 & \verb?...010101? &
\cZpLC
 & \verb?...010101? \\
\end{tabular}}

\medskip

\ZpL is also well suited for working with coefficients 
with unbalanced precision.

\smallskip

{\noindent \small
\begin{tabular}{rl}
\cIn
 & \verb?x = ?{\color{function}\verb?random_element?}\verb?(?{\color{ring}\verb?Z2?}\verb?, prec=10)? \\
 & \verb?y = ?{\color{function}\verb?random_element?}\verb?(?{\color{ring}\verb?Z2?}\verb?, prec=5)? \\
\cIn
 & \verb?u, v = x+y, x-y? \\
 & \verb?u, v? \\
\cZpCR
 & \verb?(...10111, ...01111)? \\
\cZpLC
 & \verb?(...10111, ...01111)? \\
\end{tabular}}

\smallskip

\noindent
Now, let us compute $u+v$ and compare it with $2x$ (observe that
they should be equal).

\smallskip

{\noindent \small
\begin{tabular}{rl@{\hphantom{\hspace{0.8cm}}}rl}
\cIn
 & \verb?u + v? &
\cIn
 & \verb?2*x? \\
\cZpCR
 & \verb?...00110? &
\cZpCR
 & \verb?...00110100110? \\
\cZpLC
 & \verb?...00110100110? &
\cZpLC
 & \verb?...00110100110? \\
\end{tabular}}

\smallskip

\noindent
Again \ZpCR does not output the optimal precision when the computation
is split into several steps whereas \ZpL does.
Actually, these toy examples illustrate quite common situations which 
often occur during the execution of many algorithms. For this reason,
interval arithmetic often overestimates the losses of precision. 
Roughly speaking, the aim of our package is to ``fix this misfeature''.
In the next subsections, we present a bunch of examples showing the
benefit of \ZpL in various contexts.

\medskip

\noindent \textbf{SOMOS 4.}
A first example is the SOMOS-4 sequence. It is defined by the recurrence: 
$$u_{n+4} = \frac{u_{n+1} u_{n+3} + u_{n+2}^2}{u_n}$$
and is known for its high numerical instability (see 
\cite{caruso-roe-vaccon:14a}). Nevertheless, the \ZpL package saves 
precision even when using a generic unstable implementation of the SOMOS 
iteration.

\smallskip

{\noindent \small
\begin{tabular}{rl}
\cIn
 & {\color{keyword}\verb?def?}\verb? ?{\color{function}\verb?somos4?}\verb?(u0, u1, u2, u3, n):? \\
 & \verb?    a, b, c, d = u0, u1, u2, u3? \\
 & \verb?    ?{\color{keyword}\verb?for?}\verb? _ in ?{\color{keyword}\verb?range?}\verb?(4, n+1):? \\
 & \verb?        a, b, c, d = b, c, d, (b*d + c*c) / a? \\
 & \verb?    ?{\color{keyword}\verb?return?}\verb? d? \\
\cIn
 & \verb?u0 = u1 = u2 = ?{\color{ring}\verb?Z2?}\verb?(1,15); u3 = ?{\color{ring}\verb?Z2?}\verb?(3,15)? \\
 & {\color{function}\verb?somos4?}\verb?(u0, u1, u2, u3, 18)? \\
\cZpCR
 & \verb?...11? \\
\cZpLC
 & \verb?...100000000000111? \\
\cIn
 & {\color{function}\verb?somos4?}\verb?(u0, u1, u2, u3, 100)? \\
\cZpCR
 & {\color{error}\verb?PrecisionError?}\verb?: cannot divide by something?\\
 &\verb?indistinguishable from zero.? \\
\cZpLC
 & \verb?...001001001110001? \\
\end{tabular}}

%
%

\subsection{Linear algebra}
\label{ssec:linalg}

Many generic algorithms of linear algebra lead to quite important 
instability when they are used with $p$-adic numbers. In many cases, 
our package \ZpL rubs this instability without having to change the 
algorithm, nor the implementation.

\medskip

\noindent \textbf{Matrix multiplication.}
As revealed in \cite{caruso-roe-vaccon:15}, a first simple example where 
instability appears is simply matrix multiplication.
This might be surprising because no division occurs in this situation.
Observe nevertheless the difference between \ZpCR and \ZpLC.

\smallskip

{\noindent \small
\begin{tabular}{rl}
\cIn
 & \verb?MS = ?{\color{constructor}\verb?MatrixSpace?}\verb?(?{\color{ring}\verb?Z2?}\verb?,2)? \\
 & \verb?M = ?{\color{function}\verb?random_element?}\verb?(MS, prec=5)? \\
 & \verb?for _ in range(25):? \\
 & \verb?    M *= ?{\color{function}\verb?random_element?}\verb?(MS, prec=5)? \\
 & \verb?M? \\
\cZpCR
 & \verb?[0 0]? \\
 & \verb?[0 0]? \\
\cZpLC
 & \verb?[...100000000000   ...1000000000]? \\
 & \verb?[   ...010000000     ...00100000]? \\
\end{tabular}}


\smallskip

\noindent
On the aforementioned example, we notice that \ZpCR is unable to decide 
whether the product vanishes or not. Having good estimates on the 
precision is therefore very important in such situations.

\medskip

\noindent \textbf{Characteristic polynomials.}
Characteristic polynomials are notoriously hard to compute 
\cite{caruso-roe-vaccon:15, caruso-roe-vaccon:17}. We illustrate this 
with the following example (using the default algorithm of \sage for the 
computation of the characteristic polynomial, which is a division free 
algorithm in this setting) :

\smallskip

{\noindent \small
\begin{tabular}{rl}
\cIn
 & \verb?M = ?{\color{function}\verb?random_element?}\verb?(?{\color{constructor}\verb?MatrixSpace?}\verb?(?{\color{ring}\verb?Q2?}\verb?,3), prec=10)? \\
 & \verb?M.?{\color{method}\verb?determinant?}\verb?()? \\
\cZpCR
 & \verb?...010000010? \\
\cZpLC
 & \verb?...010000010? \\
\cIn
 & \verb?M.?{\color{method}\verb?charpoly?}\verb?()? \\
\cZpCR
 & \verb?...00000000000000000001*x^3 + ? \\
 & \verb?...1001011.011*x^2 + ...0111.01*x + 0? \\
\cZpLC
 & \verb?...00000000000000000001*x^3 + ? \\
 & \verb?...1001011.011*x^2 + ...11100111.01*x +? \\
 & \verb?...010000010? \\
\end{tabular}}


\smallskip

\noindent
We observe that \ZpLC can guarantee $4$ more digits on the $x$ 
coefficient. Moreover, it recovers the correct precision on the
constant coefficient (which is the determinant) whereas \ZpCR is 
confused and cannot even certify that it does not vanish.

\subsection{Commutative algebra}
\label{ssec:comalg}

Our package can be applied to computation with $p$-adic polynomials. 

\medskip

\noindent \textbf{Euclidean algorithm.}
A natural example is that of the computation of GCD, whose stability has 
been studied in \cite{caruso:2017}.
A naive implementation of the Euclidean algorithm can produce different 
behavior depending on the type of implementation of the field of 
$p$-adic coefficients.

\smallskip

{\noindent \small
\begin{tabular}{@{}rl}
\cIn
 & \verb?S.<x> = ?{\color{constructor}\verb?PolynomialRing?}\verb?(?{\color{ring}\verb?Z2?}\verb?)? \\
 & \verb?P = ?{\color{function}\verb?random_element?}\verb?(S, degree=10, prec=5)? \\
 & \verb?Q = ?{\color{function}\verb?random_element?}\verb?(S, degree=10, prec=5)? \\
 & \verb?D = x^5 + ?{\color{function}\verb?random_element?}\verb?(S, degree=4, prec=8); D? \\
\cZpCR
 & \verb?...00000000000000000001*x^5 + ...11111010*x^4 +? \\ 
 & \verb?...10000000*x^3 + ...11001111*x^2 +? \\
 & \verb?...10000110*x + ...11100010? \\
\cZpLC
 & \verb?...00000000000000000001*x^5 + ...11111010*x^4 +? \\ 
 & \verb?...10000000*x^3 + ...11001111*x^2 +? \\
 & \verb?...10000110*x + ...11100010? \\
\cIn
 & {\color{keyword}\verb?def?}\verb? ?{\color{function}\verb?euclidean?}\verb?(A,B):? \\
 & \verb?    ?{\color{keyword}\verb?while?}\verb? B != 0:? \\
 & \verb?        A, B = B, A % B? \\
 & \verb?    ?{\color{keyword}\verb?return?}\verb? A.monic()? \\
 & {\color{function}\verb?euclidean?}\verb?(D*P, D*Q)? \\
\cZpCR
 & \verb?0*x^9 + ...1*x^8 + 0*x^7 + 0*x^6 + 0*x^5 +? \\
 & \verb?0*x^4 + 0*x^3 + ...1*x^2 + ...10*x + ...10? \\
\cZpLC
 & \verb?...00000000000000000001*x^5 + ...11111010*x^4 +? \\ 
 & \verb?...10000000*x^3 + ...11001111*x^2 +? \\
 & \verb?...10000110*x + ...11100010? \\
\end{tabular}}

\smallskip

\noindent
With high probability, $P$ and $Q$ are coprime, implying that the gcd of 
$DP$ is $DQ$ is $D$.  However, we observe that \ZpCR output a quite 
different result. The point is that, in the \ZpCR case, Euclidean 
algorithm stops prematurely because the test \verb?B != 0? fails 
too early due to the lack of precision.

\medskip

\noindent \textbf{Gröbner bases.}
Our package can be applied on complex computations like that of Gröbner 
bases using generic Gröbner bases algorithms.

\smallskip

{\noindent \small
\begin{tabular}{@{}rl}
\cIn
 & \verb?R.<x,y,z> = ?{\color{constructor}\verb?PolynomialRing?}\verb?(?{\color{ring}\verb?Q2?}\verb?, order=?{\color{string}\verb?'invlex'?}\verb?)? \\
 & \verb?F = [ ?{\color{ring}\verb?Q2?}\verb?(2,10)*x + ?{\color{ring}\verb?Q2?}\verb?(1,10)*z, ? \\
 & \verb?      ?{\color{ring}\verb?Q2?}\verb?(1,10)*x^2 + ?{\color{ring}\verb?Q2?}\verb?(1,10)*y^2 - ?{\color{ring}\verb?Q2?}\verb?(2,10)*z^2,? \\
 & \verb?      ?{\color{ring}\verb?Q2?}\verb?(4,10)*y^2 + ?{\color{ring}\verb?Q2?}\verb?(1,10)*y*z + ?{\color{ring}\verb?Q2?}\verb?(8,10)*z^2 ]? \\
\cIn
 & {\color{keyword}\verb?from?}\verb? sage.rings.polynomial.toy_buchberger\?\\
 & \verb?    ?{\color{keyword}\verb?import?}\verb? ?{\color{function}\verb?buchberger_improved?} \\
 & \verb?g = ?{\color{function}\verb?buchberger_improved?}\verb?(ideal(F))? \\
 & \verb?g.?{\color{method}\verb?sort?}\verb?(); g? \\
\cZpCR
 & \verb?[x^3, x*y + ...1100010*x^2,? \\
 & \verb? y^2 + ...11001*x^2, z + ...0000000010*x]? \\
\cZpLC
 & \verb?[x^3, x*y + ...111100010*x^2,? \\
 & \verb? y^2 + ...1111111001*x^2, z + ...0000000010*x]? \\
\end{tabular}}

\smallskip

\noindent
As we can see, some loss in precision occurs in the
Buchberger algorithm and is avoided thanks to \ZpL.

\subsection{$p$-adic differential equations}
\label{ssec:diffeq}

In \cite{LV16}, the behavior of the precision when solving
$p$-adic differential equations with separation of 
variables has been studied.
The authors have investigated the gap that appears
when applying a Newton-method solver between
the theoretic loss in precision and the 
actual loss in precision for a naive implementation in \verb?Zp(p)?.
We can reach this theoretical loss in precision using \ZpL.
We use a generic {\color{function}\verb?Newton_Iteration_Solver?}\verb?(g,h,N)?
that applies \verb?N? steps of the Newton method for 
$y'=g \times h(y)$ as described in \cite{LV16}.

\smallskip

{\noindent \small
\begin{tabular}{rl}
\cIn
 & \verb?S.<t> = ?{\color{constructor}\verb?PowerSeriesRing?}\verb?(?{\color{ring}\verb?Q2?}\verb?, 16)? \\
 & \verb?h = 1 + t + t^3? \\
 & \verb?y = t + t^2 * ?{\color{function}\verb?random_element?}\verb?(S, prec=10)? \\
 & \verb?g = y.derivative() / h(y)? \\
 & \verb?u = ?{\color{function}\verb?Newton_Iteration_Solver?}\verb?(g, h, 4); u[15]? \\
\cZpCR
 & \verb?...1101? \\
\cZpLC
 & \verb?...11011101? \\
\end{tabular}}


\section{Behind the scenes}
\label{sec:implementation}

In this section, we explain how our package \ZpL works and analyze
its performance.
The main theoretical result on which our package is based is the 
ultrametric precision theory developed in \cite{caruso-roe-vaccon:14a}, which 
suggests tracking precision \emph{via} lattices and differential 
computations. For this reason, our approach is very inspired by 
automatic differentiation techniques \cite{rall:1981} and our implementation 
follows the usual operator overloading strategy.
We will introduce two versions of our package, namely \ZpLC and \ZpLF: 
this former is safer while the latter is faster.

\smallskip

\noindent \textit{Remark about the naming.}
The letter \texttt{L}, which appears in the name of the package, 
comes from ``lattices''. The letters \texttt{C} (in \ZpLC) and 
\texttt{F} (in \ZpLF) stand for ``cap'' and ``float'' respectively.

\subsection{The precision Lemma}
\label{ssec:preclemma}

In \cite{caruso-roe-vaccon:14a}, we suggest the 
use of lattices to represent the precision of elements in 
$\Qp$-vector spaces.  This approach contrasts with the
\emph{coordinate-wise method} (of \textit{e.g.}
 \verb?Zp(5)?) that is traditionally used in \sage
where the precision of an element is specified by giving the precision
of each coordinate separately and is updated after each basic
operation.

Consider a finite dimensional normed vector space $E$ defined over 
$\Qp$. We use the notation $\Vert \cdot \Vert_E$ for the norm on $E$ and 
$B^-_E(r)$ (resp. $B^{\phantom -}_E(r)$) for the open (resp. closed) 
ball of radius $r$ centered at the origin. A \emph{lattice} $L \subset 
E$ is a sub-$\Zp$-module which generates $E$ over $\Qp$. Because of 
ultrametricity, the balls $B^{\phantom -}_E(r)$ and $B^-_E(r)$ are 
examples of lattices. Lattices can be thought of as special 
neighborhoods of $0$, and therefore are good candidates to model 
precision data. Moreover, as revealed in \cite{caruso-roe-vaccon:14a}, 
they behave quite well under (strictly) differentiable maps:

\begin{prop}
\label{prop:precision}
Let $E$ and $F$ be two finite dimensional normed vector spaces over $\Qp$ 
and $f : U \rightarrow F$ be a function defined on an open subset $U$ of 
$E$. We assume that $f$ is differentiable at some point $v_0 \in U$ and 
that the differential $df_{v_0}$ is surjective.
Then, for all $\rho \in (0, 1]$, there exists a positive real
number $\delta$ such that, for all $r \in (0, \delta)$, any lattice
$H$ such that $B^-_E(\rho r) \subset H \subset B^{\phantom -}_E(r)$ 
satisfies:
\begin{equation}
\label{eq:firstorder}
f(v_0 + H) = f(v_0) + df_{v_0} (H).
\end{equation}
\end{prop}

This proposition enables the \emph{lattice method} of tracking precision,
where the precision of the input is specified as a lattice $H$ and precision
is tracked via differentials of the steps within a given algorithm.
The equality sign in Eq.~\eqref{eq:firstorder} shows that this method
yields the optimum possible precision. 
We refer to \cite[\S 4.1]{caruso-roe-vaccon:14a} for a more complete 
exposition.

%

\subsection{Tracking precision}
\label{ssec:trackprec}

We now explain in more details the internal mechanisms \ZpLC and \ZpLF 
use for tracking precision. 

In what follows, it will be convenient to use a notion of discrete time 
represented by the letter $t$. Rigorously, it is defined as follows: 
$t=0$ when the $p$-adic ring $\ZpLC(\:\cdots)$ or $\ZpLF(\:\cdots)$ is 
created and increases by $1$ each time a variable is created, 
deleted\footnote{The deletion can be explicit (through a call to the 
\texttt{del} operator) or implicit (handled by the garbage collector).} 
or updated.

Let $\calV_t$ be the set of alive variables at time $t$. Set $E_t = 
\Qp^{\calV_t}$; it is a finite dimensional vector space over $\Qp$ which 
should be thought of as the set of all possible values that can be taken by 
the variables in $\calV_t$. For $\ttv \in \calV_t$, let $e_\ttv \in
E_t$ be the vector whose coordinates all vanish except 
at position $\ttv$ which takes the value $1$. The family 
$(e_\ttv)_{\ttv \in \calV_t}$ is obviously a basis of $E_t$; we will 
refer to it as the \emph{canonical basis}.

\subsubsection{The case of \ZpLC.}

Following Proposition~\ref{prop:precision},
the package \ZpLC follows the precision by keeping track of a lattice
$H_t$ in $E_t$, which is a global object whose purpose is to model the 
precision on all the variables in $\calV_t$ all together.
Concretely, this lattice is represented by a matrix $M_t$ in row-echelon 
form whose rows form a set of generators.
Below, we explain how the matrices $M_t$ are updated each time
$t$ increases.

\smallskip

\noindent \textit{Creating a variable.}
This happens when we encounter an instruction having one of the
two following forms:

\medskip

\noindent \hspace{5mm} \makebox[2.5cm]{[Computation]\hfill\null}
\verb?w = ?$f$\verb?(v_1,? \ldots\verb?, v_n)?

\smallskip

\noindent \hspace{5mm} \makebox[2.5cm]{[New value]\hfill\null}
\verb?w = ?\verb?R(?\text{value}\verb?, ?\text{prec}\verb?)?

\medskip

\noindent
In both cases, $\ttw$ is the newly created variable. 
The $\ttv_i$'s stand for already defined variables and $f$ is some 
$n$-ary builtin function (in most cases it is just addition, 
subtraction, multiplication or division). On the contrary, the terms 
``value'' and ``prec'' refer to user-specified constants or integral 
values which was computed earlier.

\smallskip

Let us first examine the first construction [Computation].
With our conventions, if $t$ is the time just before the execution of 
the instruction we are interested in, the $\ttv_i$'s lie in $\calV_t$ 
while $\ttw$ does not. Moreover $\calV_{t+1} = \calV_t \sqcup \{\ttw\}$, 
so that $E_{t+1} = E_t \oplus \Qp e_\ttw$.
The mapping taking the values of variables at time $t$ to that at 
time $t{+}1$ is:
$$\begin{array}{rcl}
F: \quad E_t & \longrightarrow & E_{t+1} \smallskip \\
\underline x & \mapsto & \underline x \:\oplus f(x_1, \ldots, x_n)
\end{array}$$
where $x_i$ is the $\ttv_i$-th coordinate of the vector $\underline x$.
The Jacobian matrix of $F$ at $\underline x$ is easily computed; it 
is the block matrix $J_{\underline x}(F) = \big( \begin{matrix} I & L
\end{matrix} \big)$ where $I$ is the identity matrix of size $\text{Card}
\:\calV_t$ and $L$ is the column vector whose $\ttv$-th entry is
$\frac{\partial f}{\partial \ttv} (\underline x)$ if $\ttv$ is one of 
the $\ttv_i$'s and $0$ otherwise.
Therefore, the image of $H_t$ under $dF_{\underline x}$ is represented 
by the matrix $J_{\underline x}(F) \cdot M_t = \big( \begin{matrix} M_t & C
\end{matrix} \big)$ where $C$ is the column vector:
\begin{equation}
\label{eq:defC}
C = \sum_{i=1}^n \frac{\partial f}{\partial{\ttv_i}} (\underline x)
\cdot C_i
\end{equation}
where $C_i$ is the column vector of $M_t$ corresponding to the
variable $\ttv_i$.
Observe that the matrix $J_{\underline x}(F) \cdot M_t$ is no longer a 
square matrix; it has one extra column. This reflects the fact that 
$\dim E_{t+1} = \dim E_t + 1$. Rephrasing this in a different language, 
the image of $H_t$ under $dF_{\underline x}$ is no longer a lattice in 
$E_{t+1}$ but is included in an hyperplane.

The package \ZpLC tackles this issue by introducing a cap: we do not 
work with $dF_{\underline x}(H_t)$ but instead define the lattice
$H_{t+1} = dF_{\underline x}(H_t) \oplus p^{N_{t+1}} \Zp e_\ttw$
where $N_{t+1}$ is an integer, the so-called \emph{cap}.
Alternatively, one may introduce the map:
\begin{equation}
\label{eq:tildeF}
\begin{array}{rcl}
\tilde F: \quad E_t \oplus \Qp & \longrightarrow & E_{t+1} \smallskip \\
\underline x \oplus c & \mapsto & \underline x \:\oplus 
\big(f(x_1, \ldots, x_n) + c\big).
\end{array}
\end{equation}
The lattice $H_{t+1}$ is then the image of $H_t \oplus p^{N_{t+1}}
\Zp$ under $d\tilde F_{(\underline x, \star)}$ for any value of $\star$.
The choice of the 
cap is of course a sensitive question. \ZpLC proceeds as follows. When a 
ring is created, it comes with two constants (which can be specified by 
the user): a relative cap \textsc{relcap} and an absolute cap 
\textsc{abscap}. With these predefined values, the chosen cap is:
\begin{equation}
\label{eq:cap}
N_{t+1} = 
\min\big(\textsc{abscap},\, \textsc{relcap} + v_p(y)\big)
\end{equation}
with $y = f(x_1, \ldots, x_n)$.
In concrete terms, the lattice $H_{t+1}$ is represented by the block 
matrix:
$$\left(\begin{matrix}
M_t & C \smallskip \\ 0 & p^{N_{t+1}}
\end{matrix}\right).$$
Performing row operations, we see then the entries of $C$ can be
reduced modulo $p^{N_{t+1}}$ without changing the lattice. In order
to optimize the size of the objects, we perform this reduction and
define $M_{t+1}$ by:
$$M_{t+1} = \left(\begin{matrix}
M_t & C \text{ mod } p^{N_{t+1}} \smallskip \\ 0 & p^{N_{t+1}}
\end{matrix}\right).$$
We observe in particular that $M_{t+1}$ is still in row-echelon form.

Finally, we need to explain which value is set to the newly created 
variable $\ttw$. We observe that it cannot be exactly $f(x_1, \ldots, 
x_n)$ because the latter is \emph{a priori} a $p$-adic number which 
cannot be computed exactly. For this reason, we have to truncate it at 
some finite precision. Again we choose the precision $O(p^{N_{t+1}})$,
\emph{i.e.} we define $x_\ttw$ as $f(x_1, \ldots, x_n) \text{ mod } 
p^{N_{t+1}}$. 
The congruence
$\bar x \oplus f(x_1, \ldots, x_n) \equiv
\bar x \oplus x_\ttw \pmod{H_{t+1}}$
(which holds thanks to the extra generator we have added) justifies
this choice.

\smallskip

The second construction
``\verb?w = ?\verb?R(?\text{value}\verb?, ?\text{prec}\verb?)?''
is easier to handle since, roughly speaking, it corresponds to the
case $n = 0$.
In this situation, keeping in mind the cap, the lattice $H_{t+1}$ is 
defined by
$H_{t+1} = H_t + p^{\min(\text{prec}, N_{t+1})} \Zp e_\ttw$
for the cap
$N_{t+1} = \min\big(\textsc{abscap},\, \textsc{relcap} + 
v_p(\text{value})\big)$.
The corresponding matrix $M_{t+1}$ is then given by:
$$M_{t+1} = \left(\begin{matrix}
M_t & 0 \smallskip \\ 0 & p^{\min(\text{prec}, N_{t+1})}
\end{matrix}\right).$$

\smallskip

\noindent \textit{Deleting a variable.}
Let us now examine the case where a variable $\ttw$ is deleted (or 
collected by the garbage collector). Just after the deletion, at
time $t{+}1$, we then have $\calV_{t+1} = \calV_t \backslash \{\ttw\}$.
Thus $E_t = E_{t+1} \oplus \Qp e_\ttw$. Moreover, the deletion of $\ttw$
is modeled by the canonical projection $f : E_t \to E_{t+1}$. Since $f$
is linear, it is its own differential (at each point) and we set 
$H_{t+1} = f(H_t)$.
A matrix representing $H_{t+1}$ is deduced from $M_t$ by erasing the 
column corresponding to $\ttw$. However the matrix we get this way is 
no longer in row-echelon form. We then need to re-echelonize it.

More precisely, the obtained matrix has this shape:
$$\left(\raisebox{-0.5\height}
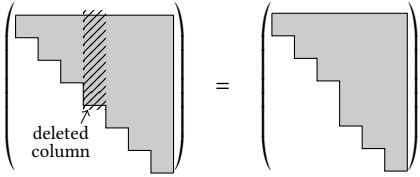
{\begin{tikzpicture}[scale=0.3]
\fill[black!20] (0,0)--(0,-1)--(1,-1)--(1,-2)--(2,-2)--(2,-3)
    --(3,-3)--(3,-4)--(4,-4)--(4,-5)--(5,-5)--(5,-6)--(6,-6)
    --(6,-7)--(7,-7)--(7,0)--cycle;
\draw (0,0)--(0,-1)--(1,-1)--(1,-2)--(2,-2)--(2,-3)
    --(3,-3)--(3,-4)--(4,-4)--(4,-5)--(5,-5)--(5,-6)--(6,-6)
    --(6,-7)--(7,-7);
\draw[thin] (7,-7)--(7,0)--(0,0);
\begin{scope}
\clip (3,0.2) rectangle (4,-4.2);
\draw[pattern=north east lines] (2.5,0.5) rectangle (4.5,-4.5);
\end{scope}
\draw[->] (2.9,-4.8)--(3.2,-4.2);
\node[scale=0.8] at (2,-5.2) { deleted };
\node[scale=0.8] at (2,-6) { column };
\end{tikzpicture}}\right)
\quad = \quad
\left(\raisebox{-0.5\height}{\begin{tikzpicture}[scale=0.3]
\fill[black!20] (0,0)--(0,-1)--(1,-1)--(1,-2)--(2,-2)--(2,-3)
    --(3,-3)--(3,-5)--(4,-5)--(4,-6)--(5,-6)--(5,-7)--(6,-7)--(6,0)--cycle;
\draw (0,0)--(0,-1)--(1,-1)--(1,-2)--(2,-2)--(2,-3)
    --(3,-3)--(3,-5)--(4,-5)--(4,-6)--(5,-6)--(5,-7)--(6,-7);
\draw[thin] (6,-7)--(6,0)--(0,0);
\end{tikzpicture}}\right)$$
where a cell is colored when it can contain a non-vanishing entry.
The top part of the matrix is then already echelonized, so that we
only have to re-echelonize the bottom right corner whose size is
the distance from the column corresponding to the erased variable
to the end. 
Thanks to the particular shape of the matrix, the echelonization can be 
performed efficiently: we combine the first rows (of the bottom right 
part) in order to clear the first unwanted nonzero entry and then 
proceed recursively.

\medskip

\noindent \textit{Updating a variable.}
Just like for creation, this happens when the program reaches an
affectation ``\verb?w = ...?'' where the variable \ttw is already 
defined. This situation reduces to the creation of the temporary
variable (the value of the right-hand-size), the deletion of the
old variable \ttw and a renaming. It can then be handled using the
methods discussed previously.

\subsubsection{The case of \ZpLF.}

The way the package \ZpLF tracks precision is based on similar 
techniques but differs from \ZpLC in that it does not introduce a cap 
but instead allows $H_t$ to be a sub-$\Zp$-module of $E_t$ of any 
codimension.
This point of view is nice because it implies smaller objects and 
consequently leads to faster algorithms. However, it has a huge 
drawback; indeed, unlike lattices, submodules of $E_t$ of arbitrary 
codimensions are \emph{not} exact objects, in the sense that they cannot 
be represented by integral matrices in full generality. Consequently,
they cannot be encoded on a computer.
We work around this drawback by replacing everywhere exact $p$-adic 
numbers by floating point $p$-adic numbers (at some given precision) 
\cite{course-padic}.

The fact that the lattice $H_t$ can now have arbitrary codimension
translates to the fact the matrix $M_t$ can be rectangular.
Precisely, we will maintain matrices $M_t$ of the shape:
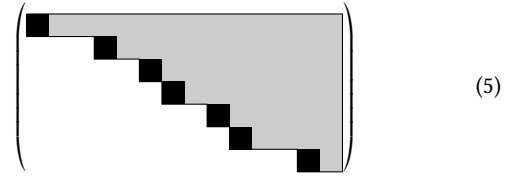
\begin{equation}
\label{eq:shapeM-ZpLF}
\left(\raisebox{-0.5\height}{\begin{tikzpicture}[scale=0.3]
\fill[black!20] (0,0)--(0,-1)--(3,-1)--(3,-2)--(5,-2)--(5,-3)
    --(6,-3)--(6,-4)--(8,-4)--(8,-5)--(9,-5)--(9,-6)--(12,-6)
    --(12,-7)--(14,-7)--(14,0)--cycle;
\draw (0,0)--(0,-1)--(3,-1)--(3,-2)--(5,-2)--(5,-3)
    --(6,-3)--(6,-4)--(8,-4)--(8,-5)--(9,-5)--(9,-6)--(12,-6)
    --(12,-7)--(14,-7);
\fill (0,0) rectangle (1,-1);
\fill (3,-1) rectangle (4,-2);
\fill (5,-2) rectangle (6,-3);
\fill (6,-3) rectangle (7,-4);
\fill (8,-4) rectangle (9,-5);
\fill (9,-5) rectangle (10,-6);
\fill (12,-6) rectangle (13,-7);
\draw[thin] (14,-7)--(14,0)--(0,0);
\end{tikzpicture}}\right)
\end{equation}
where only the colored cells may contain a nonzero value and the 
black cells ---the so-called \emph{pivots}--- do not vanish. A
variable whose corresponding column contains a pivot will be called 
a \emph{pivot variable at time $t$}.


\medskip

\noindent \textit{Creating a variable.}
We assume first that the newly created variable is defined through a
statement of the form:
``\verb?w = ?$f$\verb?(v_1,? \ldots\verb?, v_n)?''.
As already explained in the case of \ZpLC, this code is 
modeled by the mathematical mapping:
$$\begin{array}{rcl}
F: \quad E_t & \longrightarrow & E_{t+1} \smallskip \\
\underline x & \mapsto & \underline x \:\oplus f(x_1, \ldots, x_n).
\end{array}$$
Here $\underline x$ represents the state of memory at time $t$, 
and $x_i$ is the coordinate of $\underline x$ corresponding to the
variable $\ttv_i$.

In the \ZpLF framework, $H_{t+1}$ is defined as the image of $H_t$
under the differential $dF_{\underline x}$. Accordingly, the matrix
$M_{t+1}$ is defined as $M_{t+1} = \left(\begin{matrix}
M_t & C \end{matrix}\right)$
where $C$ is the column vector defined by Eq.~\eqref{eq:defC}.
However, we insist on the fact that all the computations now take
place in the ``ring'' of floating point $p$-adic numbers. Therefore,
we cannot guarantee that the rows of $M_{t+1}$ generate $H_{t+1}$.
Nonetheless, they generate a module which is expected to be close
to $H_{t+1}$.

\smallskip

If \ttw is created by the code
``\verb?w = ?\verb?R(?\text{value}\verb?, ?\text{prec}\verb?)?'',
we define $H_{t+1} = H_t \oplus p^{\text{prec}} \Zp e_\ttw$ and consequently:
$$M_{t+1} = \left(\begin{matrix}
M_t & 0 \smallskip \\ 0 & p^{\text{prec}}
\end{matrix}\right)$$
If $\text{prec}$ is $+\infty$ (or, equivalently, not specified), we
agree that $H_{t+1} = H_t$ and $M_{t+1} = (\begin{matrix} M_t & 0 
\end{matrix})$.

\medskip

\noindent \textit{Deleting a variable.} 
As for \ZpLC, the matrix operation implied by the deletion of the 
variable \ttw is the deletion of the corresponding column of $M_t$. If 
\ttw is not a pivot variable at time $t$, the matrix $M_t$ keeps the 
form \eqref{eq:shapeM-ZpLF} after erasure; therefore no more treatment 
is needed in this case.

Otherwise, we re-echelonize the matrix as follows. After the deletion 
of the column $C_\ttw$, we examine the first column $C$ which was 
located on the right of $C_\ttw$. Two situations may occur (depending
on the fact that $C$ was or was not a pivot column):

\noindent \hfill
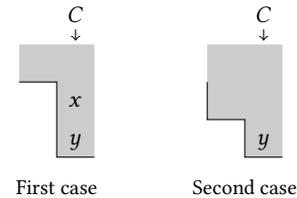
\begin{tikzpicture}[scale=0.5]
\fill[black!20] (0,0)--(1,0)--(1,-2)--(2,-2)--(2,1)--(0,1)--(0,0);
\draw (0,0)--(1,0)--(1,-2)--(2,-2);
\node at (1.5,1.8) { $C$ };
\draw[->] (1.5,1.4)--(1.5,1.1);
\node at (1.5,-0.5) { \ph $x$ };
\node at (1.5,-1.5) { \ph $y$ };
\node[scale=0.9] at (1,-2.8) { First case };
\begin{scope}[xshift=5cm]
\fill[black!20] (0,0)--(0,-1)--(1,-1)--(1,-2)--(2,-2)--(2,1)--(0,1)--(0,0);
\draw (0,0)--(0,-1)--(1,-1)--(1,-2)--(2,-2);
\node at (1.5,1.8) { $C$ };
\draw[->] (1.5,1.4)--(1.5,1.1);
\node at (1.5,-1.5) { \ph $y$ };
\node[scale=0.9] at (1,-2.8) { Second case };
\end{scope}
\end{tikzpicture}
\hfill\null

\noindent

In the first case, we perform row operations in order to replace the 
pair $(x,y)$ by $(d, 0)$ where $d$ is an element of valuation 
$\min(v_p(x), v_p(y))$. Observe that $y$ is necessarily nonzero in this case, 
so that $d$ does not vanish as well. After this operation, we move to
the next column and repeat the same process.

The second case is divided into two subcases. First, if $y$ does not 
vanish, it can serve as a pivot and the obtained matrix has the desired 
shape. When this occurs, the echelonization stops. On the contrary, if 
$y = 0$, we just untint the corresponding cell and move to the next 
column without modifying the matrix.

\subsection{Visualizing the precision}
\label{ssec:viewprec}

Our package implements several methods giving access to the precision 
structure. In the subsection, we present and discuss the most relevant 
features in this direction.

\medskip

\noindent \textbf{Absolute precision of one element.}
This is the simplest accessible precision datum.
It is encapsulated in the notation when an element is printed. For
example, the (partial) session:

\smallskip

{\noindent \small
\noindent
\begin{tabular}{rl}
\cIn   & \verb?v = ?{\color{ring}\verb?Z2?}\verb?(173,10); v? \\
\cZpLC & \verb?...0010101101?
\end{tabular}}

\smallskip

\noindent
indicates that the absolute precision on \ttv is $10$ since
exactly $10$ digits are printed.
The method {\color{method}\verb?precision_absolute?} provides a more easy-to-use
access to the absolute precision.

\smallskip

{\noindent \small
\noindent
\begin{tabular}{rl}
\cIn   & \verb?v.?{\color{method}\verb?precision_absolute?}\verb?()? \\
\cZpLC & \verb?10?
\end{tabular}}

\smallskip

\noindent
Both \ZpLC and \ZpLF compute the absolute precision of \ttv (at time 
$t$) as the smallest valuation of an entry of the column of $M_t$ 
corresponding to the variable \ttv. Alternatively, it is
the unique integer $N$ for which $\pi_\ttv(H_t) = p^N \Zp$ where 
$\pi_\ttv : E_t \to \Qp$ takes a vector to its \ttv-coordinate.
This definition of the absolute precision sounds revelant because, if we 
believe that the submodule $H_t \subset E_t$ is supposed to encode the 
precision on the variables in $\calV_t$, 
Proposition~\ref{prop:precision} applied with the mapping $\pi_\ttv$ 
indicates that a good candidate for the precision on $e_\ttv$ is 
$\pi_\ttv(H_t)$, that is $p^N \Zp$.

\medskip

\noindent \textit{About correctness.}
We emphasize that the absolute precision computed this way is \emph{not} 
proved, either for \ZpLF or \ZpLC. However, in the case of \ZpLC, one
can be slightly more precise. Let $\calU_t$ be the vector space of 
user-defined variables before time $t$ and $U_t$ be the lattice
modeling the precision on them. The pair $(\calU_t, U_t)$ is defined 
inductively as follows: we set $\calU_0 = U_0 = 0$ and 
$\calU_{t+1} = \calU_t \oplus \Qp e_\ttw$,
$U_{t+1} = U_t \oplus p^{\text{prec}} \Zp e_\ttw$
when a new variable \ttw is created by
``\verb?w = R(?value\verb?, ?prec\verb?)?''; otherwise, we put 
$\calU_{t+1} = \calU_t$ and $U_{t+1} = U_t$. Moreover the values
entered by the user defines a vector (with integral coordinates)
$\underline u_t \in \calU_t$.

Similarly, in order to model the caps, we define a pair $(\calK_t, 
K_t)$ by the recurrence $\calK_{t+1} = \calK_t \oplus \Qp e_\ttw$,
$K_{t+1} = K_t \oplus p^{N_{t+1}} \Zp e_\ttw$
each time a new variable \ttw is created. 
Here, the exponent $N_{t+1}$ is the cap defined by Eq.~\eqref{eq:cap}.
In case of deletion, we put $\calK_{t+1} = \calK_t$ and $K_{t+1} = K_t$.

Taking the compositum of all the functions $\tilde F$ (\emph{cf} 
Eq.~\eqref{eq:tildeF}) from time $0$ to $t$, we find that the execution 
of the session until time $t$ is modeled by a mathematical function
$\Phi_t : \calU_t \oplus \calK_t \to E_t$.
From the design of \ZpLC, we deduce further that there exists a 
vector $\underline k_t \in K_t$ such that:
$$\Phi_t(\underline u_t \oplus \underline k_t) = \underline x_t 
\quad \text{and} \quad
d\Phi_t (U_t \oplus K_t) = H_t$$
where the differential of $\Phi_t$ is taken at the point
$\underline u_t \oplus \underline k_t$. 
Set $\Phi_{t,\ttv} = \pi_\ttv \circ \Phi_t$; it maps 
$\underline u_t \oplus \underline k_t$ to the \ttv-coordinate 
$x_{t,\ttv}$ of $\underline x_t$ and satisfies
$d\Phi_{t,\ttv} (U_t \oplus K_t) = \pi_\ttv(H_t) = p^N \Zp$
where $N$ is the value returned by {\color{method}\verb?precision_absolute?}. Thus, 
as soon as the assumptions of Proposition \ref{prop:precision} are
fulfilled, we derive
$\Phi_{t,\ttv} \big((\underline u_t + U_t) \oplus (\underline k_t + K_t)\big) = 
x_{t,\ttv} + p^N \Zp$.
Noting that $k_t \in K_t$, we finally get:
\begin{equation}
\label{eq:absprec}
\Phi_{t,\ttv} (\underline u_t + U_t) \subset
\Phi_{t,\ttv} \big((\underline u_t + U_t) \oplus K_t\big) = 
x_{t,\ttv} + p^N \Zp.
\end{equation}
The latter inclusion means that the computed value $x_{t,\ttv}$ is
accurate at precision $O(p^N)$, \emph{i.e.} that the output absolute
precision is correct. 

Unfortunately, checking automatically the assumptions of Proposition 
\ref{prop:precision} in full generality seems to be difficult, though it 
can be done by hand for many particular examples
\cite{caruso-roe-vaccon:14a, caruso:2017, LV16}.

\begin{rmk}
Assuming that Proposition \ref{prop:precision} applies, 
the absolute precision computed as above is optimal if and only if
the inclusion of \eqref{eq:absprec} is an equality. Applying again
Proposition \ref{prop:precision}
with the restricted mapping $\Phi_{t,\ttv} : \calU_t \to \Qp$ and the
lattice $U_t$, we find that this happens if and only if
$d\Phi_{t,\ttv}(U_t) = p^N \Zp$.

Unfortunately, the latter condition cannot be checked on the matrix 
$M_t$ (because of reductions). 
However it is possible (and easy) to check whether the weaker condition 
$d\Phi_{t,\ttv}(K_t) \subsetneq p^N \Zp$. This checking is achieved by 
the method
{\color{method}\verb?is_precision_capped?} (provided by our package) which returns true 
if $d\Phi_{t,\ttv}(K_t) = p^N \Zp$. As a consequence, when this method
answers \textsc{false}, the absolute precision computed by the software 
is likely optimal.
\end{rmk}

\noindent \textbf{Precision on a subset of elements.}
Our package implements the method {\color{method}\verb?precision_lattice?} through which 
we can have access to the joint precision on a set of variables: it 
outputs a matrix (in echelon form) whose rows generate a lattice 
representing the precision on the subset of given variables.

When the variables are ``independent'', the precision lattice is
split and the method {\color{method}\verb?precision_lattice?} outputs a diagonal
matrix:

\smallskip

{\noindent \small
\begin{tabular}{rl}
\cIn   & \verb?x = ?{\color{ring}\verb?Z2?}\verb?(987,10); y = ?{\color{ring}\verb?Z2?}\verb?(21,5)? \\
\cIn   & {\color{comment}\verb?#? \textit{We first retrieve the precision object}} \\
       & \verb?L = ?{\color{ring}\verb?Z2?}\verb?.?{\color{method}\verb?precision?}\verb?()?\\
\cIn   & \verb?L.?{\color{method}\verb?precision_lattice?}\verb?([x,y])? \\
\cZpLC & \verb?[1024    0]? \\
       & \verb?[   0   32]?
\end{tabular}}

\smallskip

\noindent
However, after some computations, the precision matrix evolves and 
does not remain diagonal in general (though it is always triangular
because it is displayed in row-echelon form):

\smallskip

{\noindent \small
\begin{tabular}{rl}
\cIn   & \verb?u, v = x+y, x-y? \\
       & \verb?L.?{\color{method}\verb?precision_lattice?}\verb?([u,v])? \\
\cZpLC & \verb?[  32 2016]? \\
       & \verb?[   0 2048]?
\end{tabular}}

\smallskip

\noindent
The fact that the precision matrix is no longer diagonal indicates
that some well-chosen linear combinations of $u$ and $v$ are known
with more digits than $u$ and $v$ themselves. In this particular
example, the sum $u+v$ is known at precision $O(2^{11})$ while the
(optimal) precision on $u$ and $v$ separately is only $O(2^5)$.

\smallskip

{\noindent \small
\begin{tabular}{rl}
\cIn   & \verb?u, v? \\
\cZpLC & \verb?(...10000, ...00110)? \\
\cIn   & \verb?u + v? \\
\cZpLC & \verb?...11110110110? \\
\end{tabular}}

\medskip

\noindent \textit{Diffused digits of precision.}
The phenomenon observed above is formalized by the notion of diffused
digits of precision introduced in \cite{caruso-roe-vaccon:15}.
We recall briefly its definition.

\begin{deftn}
\label{def:diffused}
Let $E$ be a $\Qp$-vector space endowed with a distinguished basis 
$(e_1, \ldots, e_n)$ and write $\pi_i : E \to \Qp e_i$ for the 
projections.
Let $H \subset E$ be a lattice. The number of \emph{diffused digits of 
precision} of $H$ is the length of $H_0/H$ where $H_0 = \pi_1(H) \oplus 
\cdots \oplus \pi_n(H)$.
\end{deftn}

If $H$ represents the actual precision on some object, then
$H_0$ is the smallest diagonal lattice containing $H$. It then
corresponds to the maximal \emph{coordinate-wise} precision we can 
reach on the set of $n$ variables corresponding to the basis $(e_1,
\ldots, e_n)$.

The method {\color{method}\verb?number_of_diffused_digits?} computes the number of
diffused digits of precision on a set of variables. Observe:

\smallskip

{\noindent \small
\begin{tabular}{rl}
\cIn   & \verb?L.?{\color{method}\verb?number_of_diffused_digits?}\verb?([x,y])? \\
\cZpLC & \verb?0? \\
\cIn   & \verb?L.?{\color{method}\verb?number_of_diffused_digits?}\verb?([u,v])? \\
\cZpLC & \verb?6? \\
\end{tabular}}

\smallskip

\noindent
For the last example, we recall that the relevant precision lattice 
$H$ is generated by the $2 \times 2$ matrix:
$$\left(\begin{matrix} 2^5 & 2016 \\ 0 & 2^{11} \end{matrix}\right).$$
The minimal diagonal suplattice $H_0$ of $H$ is generated by the scalar matrix 
$2^5 \cdot \text{I}_2$ and contains $H$ with index $2^6$ in it. This is 
where the $6$ digits of precision come from.
There are easily visible here: the sum $u+v$ is known with $11$ digits,
that is exactly $6$ more digits than the summands $u$ and $v$.

\subsection{Complexity}
\label{ssec:complexity}

We now discuss the cost of the above operations.
In what follows, we shall count operations in $\Qp$. Although $\Qp$
is an inexact field, our model of complexity makes sense because the
size of the $p$-adic numbers we manipulate will all have roughly the
same size: for \ZpLF, it is the precision we use for floating point 
arithmetic while, for \ZpLC, it is the absolute cap which
was fixed at the beginning.

It is convenient to introduce a total order on $\calV_t$: for $\ttv, 
\ttw \in \calV_t$, we say that $\ttv <_t \ttw$ if \ttv was created 
before~\ttw. By construction, the columns of the matrix $M_t$ are 
ordered with respect to $<_t$. We denote by $r_t$ (resp. $c_t$) the
number of rows (resp. columns) of $M_t$. By construction $r_t$ is
also the cardinality of $\calV_t$. We have $c_t \leq r_t$ and the 
equality always holds in the \ZpLC case.

For $\ttv \in \calV_t$, we define the \emph{index} of \ttv, denoted by 
$\ind_t(\ttv)$ as the number of elements of $\calV_t$ which are not 
greater than \ttv. If we sort the elements of $\calV_t$ by increasing 
order, \ttv then appears in $\ind_t(\ttv)$-th position.
We also define the \emph{co-index} of \ttv by
$\coind_t(\ttv) = r_t - \ind_t(\ttv)$.

Similarly, for any variable $\ttv \in \calV_t$, we define the 
\emph{height} (resp. the \emph{co-height}) of \ttv at time $t$ as the 
number of pivot variables \ttw such that $\ttw \leq_t \ttv$ (resp. 
$\ttw >_t \ttv$). We denote it by $\hgt_t(\ttv)$ (resp. by 
$\cohgt_t(\ttv)$). Clearly $\hgt_t(\ttv) + \cohgt_t(\ttv) = c_t$.
The height of \ttv is the height of the significant 
part of the column of $M_t$ which corresponds to \ttv. In the case of
\ZpLC, all variables are pivot variables and thus $\hgt_t(\ttv) = 
\ind_t(\ttv)$ and $\cohgt_t(\ttv) = \coind_t(\ttv)$ for all \ttv.

\medskip

\noindent \textbf{Creating a variable.}
With the notations of \S \ref{ssec:trackprec}, it is obvious that 
creating a new variable \ttw requires:
$$O\bigg(\sum_{i=1}^n \hgt_i(\ttv_i)\bigg) \subset O(n \: c_t)$$ 
operations in $\Qp$. Here, we recall that $n$ is the arity of the 
operation defining \ttw. In most cases it is $2$; thus the above
complexity reduces to $O(c_t)$.

In the \ZpLF context, $c_t$ counts the number of user-defined variables. 
It is then expected to be constant (roughly equal to the size of the 
input) while running a given algorithm.

On the contrary, in the \ZpLC context, $c_t$ counts the number of
variables which are alive at time $t$. It is no longer expected to
be constant but evolves continuously when the algorithm runs.
\begin{figure}

\noindent\hfill%
\begin{tabular}{|l|r|r|r|r|r|}
\hline
\textbf{Dimension} &
 2 &   5 &   10 &    20 &      50 \\
\hline
\textbf{Total} &
35 & 424 & 5\:539 & 83\:369 & 3\:170\:657 \\
\hline
\textbf{Simult.} &
17 &  65 &  225 &   845 &    5\:101 \\
\hline
\end{tabular}%
\hfill\null

\vspace{1mm}

\noindent\hfill%
{\small Computation of characteristic polynomial}%
\hfill\null

\bigskip


\noindent\hfill%
\begin{tabular}{|l|r|r|r|r|r|r|}
\hline
\textbf{Degree} &
 2 &   5 &   10 &   20 &   50 &   100 \\
\hline
\textbf{Total} &
54 & 130 &  332 & 1\:036 & 4\:110 & 10\:578 \\ 
\hline
\textbf{Simult.} &
18 &  31 &   56 &  106 &  256 &   507 \\
\hline
\end{tabular}%
\hfill\null

\vspace{1mm}

\noindent\hfill%
{\small Naive Euclidean algorithm}%
\hfill\null

\caption{Numbers of involved variables}
\label{fig:variables}
\end{figure}
The tables of Figure~\ref{fig:variables} show the total number of 
created variables (which reflects the complexity) together with the 
maximum number of variables alive at the same time (which reflects 
the memory occupation) while executing two basic computations.
The first one is the computation of the characteristic polynomial of a 
square matrix by the default algorithm used by \sage for $p$-adic fields 
(which is a division-free algorithm of quartic complexity) while the
second one is the computation of the gcd of two polynomials using a
naive Euclidean algorithm (of quadratic complexity). We can observe 
that, for both of them, the memory usage is roughly equal to the 
square root of the complexity.

\medskip

\noindent \textbf{Deleting a variable.}
The deletion of the variable \ttw induces the deletion of the 
corresponding column of $M_t$, possibly followed by a partial 
row-echelonization. In terms of algebraic complexity, the deletion is 
free. The cost of the echelonization is within
$O\big(\coind_t(\ttw) \cdot \cohgt_t(\ttw) \big)$
operations in $\Qp$.

In the \ZpLF case, we expect that, most of the time, the deleted
variables were created after all initial variables were set by the
user. This means that we expect $\cohgt_t(\ttw)$ to vanish and so,
the corresponding cost to be negligible.

In the \ZpLC case, we always have $\cohgt_t(\ttw) = \coind_t(\ttw)$, so 
that the cost becomes $O\big(\coind_t(\ttw)^2\big)$. This does not look 
nice \emph{a priori}. However, the principle of temporal locality 
\cite{denning:2005}
asserts that $\coind_t(\ttw)$ tends to be small in general: destroyed 
variables are often variables that were created recently. As a basic 
example, variables which are local to a small piece of code (\emph{e.g.} 
a short function or a loop) do not survive for a long time. It turns
out that this behavior is typical in many implementations!
\begin{figure}
\noindent\hfill%
\begin{tikzpicture}[xscale=0.35,yscale=0.08]
\draw[fill=blue] (-0.45,0) rectangle (0.45,46);
\draw[fill=blue] (0.65,0) rectangle (1.45,69);
\draw[fill=blue] (1.65,0) rectangle (2.45,5);
\draw[fill=blue] (2.65,0) rectangle (3.45,7);
\draw[fill=blue] (3.65,0) rectangle (4.45,7);
\draw[fill=blue] (4.65,0) rectangle (5.45,12);
\draw[fill=blue] (5.65,0) rectangle (6.45,7);
\draw[fill=blue] (6.65,0) rectangle (7.45,5);
\draw[fill=blue] (7.65,0) rectangle (8.45,1);
\draw[fill=blue] (8.65,0) rectangle (9.45,1);
\draw (9.65,0)--(10.45,0);
\draw[fill=blue] (10.65,0) rectangle (11.45,1);
\draw[fill=blue] (11.65,0) rectangle (12.45,1);
\draw[fill=blue] (12.65,0) rectangle (13.45,2);
\draw[fill=blue] (13.65,0) rectangle (14.45,2);
\draw[fill=blue] (14.65,0) rectangle (15.45,2);
\draw[fill=blue] (15.65,0) rectangle (16.45,2);
\draw[fill=blue] (16.65,0) rectangle (17.45,1);
\draw[fill=blue] (17.65,0) rectangle (18.45,1);
\draw[fill=blue] (18.65,0) rectangle (19.45,1);
\node[scale=0.7] at (0,-2.5) { $0$ };
\node[scale=0.7] at (5,-2.5) { $5$ };
\node[scale=0.7] at (10,-2.5) { $10$ };
\node[scale=0.7] at (15,-2.5) { $15$ };
\draw (-2,0)--(-2,72);
\draw (-2.2,0)--(-1.8,0);
\node[scale=0.7,left] at (-2.1,0) { $0$ };
\draw (-2.2,10)--(-1.8,10);
\node[scale=0.7,left] at (-2.1,10) { $10$ };
\draw (-2.2,20)--(-1.8,20);
\node[scale=0.7,left] at (-2.1,20) { $20$ };
\draw (-2.2,30)--(-1.8,30);
\node[scale=0.7,left] at (-2.1,30) { $30$ };
\draw (-2.2,40)--(-1.8,40);
\node[scale=0.7,left] at (-2.1,40) { $40$ };
\draw (-2.2,50)--(-1.8,50);
\node[scale=0.7,left] at (-2.1,50) { $50$ };
\draw (-2.2,60)--(-1.8,60);
\node[scale=0.7,left] at (-2.1,60) { $60$ };
\draw (-2.2,70)--(-1.8,70);
\node[scale=0.7,left] at (-2.1,70) { $70$ };
\end{tikzpicture}%
\hfill\null

\caption{The distribution of $\coind_t(\ttw)$}
\label{fig:deletion}
\end{figure}
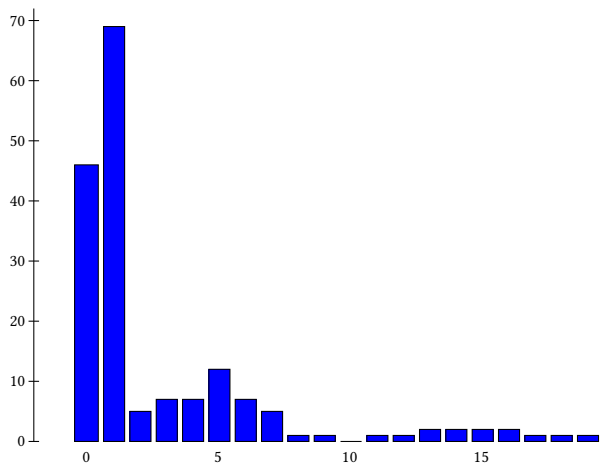
The histogram of Figure~\ref{fig:deletion} shows 
the distribution of $\coind_t(\ttw)$ while executing the Euclidean
algorithm (naive implementation) with two polynomials of degree $7$
as input. The bias is evident: most of the time $\coind_t(\ttw) \leq
1$.

\medskip

\noindent \textbf{Summary: Impact on complexity.}
We consider the case of an algorithm with the following characteristics: 
its complexity is $c$ operations in $\Qp$ (without any tracking of 
precision), its memory usage is $m$ elements of $\Qp$, its input and its 
output have size $s_\inp$ and $s_\out$ (elements of $\Qp$) respectively.

In the case of \ZpLF, creating a variable has a cost $O(s_\inp)$ 
whereas deleting a variable is free. 
Thus when executed with the \ZpLF mechanism, the complexity of our
algorithm becomes $O(s_\inp c)$.

In the \ZpLC framework, creating a variable has a cost $O(m)$. The case 
of deletion is more difficult to handle. However, by the temporal 
locality principle, it seems safe to assume that it is not the 
bottleneck (which is the case in practice). Therefore, when executed 
with the \ZpLF mechanism, the cost of our algorithm is expected to be 
roughly $O(mc)$. Going further in speculation, we might estimate the 
magnitude of $m$ as about $s + \sqrt c$ with $s = \max(s_\inp, 
s_\out)$, leading to a complexity of $O(c^{3/2} + sc)$. For 
quasi-optimal algorithms, the term $sc \simeq c^2$ dominates. However, 
as soon as the complexity is at least quadratic in $s$, the dominant 
term is $c^{3/2}$ and the impact on the complexity is then limited.

\section{Conclusion}
\label{sec:conclusion}

The package \ZpL provides powerful tools (based on automatic 
differentiation) to track precision in the $p$-adic setting. In many 
concrete situations, it greatly outperforms standard interval arithmetic,
as shown in \S \ref{sec:demo}.
The impact on complexity is controlled but nevertheless non-negligible 
(see \S \ref{ssec:complexity}). For this reason, it is unlikely that a fast
algorithm will rely \emph{directly} on the machinery proposed by \ZpL,
though it might do so for a specific part of a computation.
At least for now, bringing together rapidity and stability
still requires a substantial human contribution 
and a careful special study of all parameters.

Nevertheless, we believe that \ZpL can be extremely helpful to 
anyone designing a fast and stable $p$-adic algorithm for a couple 
of reasons. First, it provides mechanisms 
to automatically detect which steps of a given algorithm are stable and 
which ones are not. In this way, it highlights the parts of the algorithm 
on which the researcher has to concentrate their effort.
Second, recall that a classical strategy to improve stability consists 
in working internally at higher precision. Finding the internal increase
in precision that best balances efficiency and accuracy is not an easy
task in general. Understanding the number of diffused digits of 
precision gives very useful hints in this direction. For example, when there are
no diffused digits of precision then the optimal precision completely splits over the
variables and there is no need to internally increase the precision. On 
the contrary, when there are many diffused digits of precision, a large 
increment is often required.
Since \ZpL gives a direct access to the number of diffused digits of 
precision, it could be very useful to the designer who is concerned
with the balance between efficiency and accuracy.

\bibliographystyle{plain}
\bibliography{ZpL}

\end{document}